\documentclass{gtart}     
\usepackage{amssymb} 
\newtheorem{theorem}{Theorem}[section]

\newtheorem{lemma}[theorem]{Lemma}

\theoremstyle{remark}
\newtheorem{remark}[theorem]{Remark}
\newtheorem{definition}[theorem]{Definition}
\newtheorem{example}[theorem]{Example}
\begin{document}
%
\input gtoutput
\volumenumber{2}\papernumber{2}\volumeyear{1998}
\pagenumbers{11}{29}\published{19 March 1998}
\proposed{Jean-Pierre Otal}\seconded{Cameron Gordon, Walter Neumann}
\received{21 February 1997}\revised{13 March 1998}
\accepted{19 March 1998}
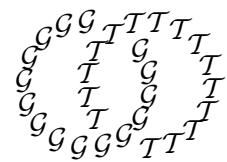
\title{The Symmetry of Intersection Numbers\\in Group Theory}

\author{Peter Scott}

\address{Mathematics Department\\University of Michigan\\Ann Arbor, MI
48109, USA} 

\email{pscott@math.lsa.umich.edu}

\begin{abstract}%
For suitable subgroups of a finitely generated group, we define the
intersection number of one subgroup with another subgroup and show that this
number is symmetric. We also give an interpretation of this number.
\end{abstract}

\primaryclass{20F32}

\secondaryclass{20E06, 20E07, 20E08, 57M07}

\keywords{Ends, amalgamated free products, trees}

\maketitlepage



If one considers two simple closed curves $L$ and $S$ on a closed orientable
surface $F,$ one can define their intersection number to be the least number
of intersection points obtainable by isotoping $L$ and $S$ transverse to
each other. (Note that the count is to be made without any signs attached to
the intersection points.) By definition, this number is symmetric, ie the
roles of $L$ and $S$ are interchangeable. This can be regarded as a
definition of the intersection number of the two infinite cyclic subgroups $%
\Lambda $ and $\Sigma $ of the fundamental group of $F$ which are carried by 
$L$ and $S.$ In this paper, we show that an analogous definition of
intersection number of subgroups of a group can be given in much greater
generality and proved to be symmetric. We also give an interpretation of
these intersection numbers.

In \cite{Rips-Sela}, Rips and Sela considered a torsion free finitely
presented group $G$ and infinite cyclic subgroups $\Lambda $ and $\Sigma $
such that $G$ splits over each. (A group $G$ splits over a subgroup $C$ if
either $G$ has a HNN\ decomposition $G=A\ast _{C},$ or $G$ has an
amalgamated free product structure $G=A\ast _{C}B,$ where $A\neq C\neq B.)$
They effectively considered the intersection number $i(\Lambda ,\Sigma )$ of 
$\Lambda $ with $\Sigma ,$ and they proved that $i(\Lambda ,\Sigma )=0$ if
and only if $i(\Sigma ,\Lambda )=0.$ Using this, they proved that $G$ has
what they call a JSJ decomposition. If $i(\Lambda ,\Sigma )$ was not zero,
it follows from their work that $G$ can be expressed as the fundamental
group of a graph of groups with some vertex group being a surface group $H$
which contains $\Lambda $ and $\Sigma .$ Now it is intuitively clear (and we
discuss it further at the end of section \ref{intersection numbers} of this
paper) that the intersection number of $\Lambda $ with $\Sigma $ is the same
whether it is measured in $G$ or in $H.$ Also the intersection numbers of $%
\Lambda $ and $\Sigma $ in $H$ are symmetric because of their topological
interpretation. So it follows at the end of all their work that the
intersection numbers of $\Lambda $ and $\Sigma $ in $G$ are also symmetric.
In 1994, Rips asked if there was a simpler proof of this symmetry which does
not depend on their proof of the JSJ splitting. The answer is positive, and
the ideas needed for the proof are all essentially contained in earlier
papers of the author. This paper is a belated response to Rips' question.
The main idea is to reduce the natural, but not clearly symmetric,
definition of intersection number to counting the intersections of suitably
chosen sets. The most general possible algebraic situation in which to
define intersection numbers seems to be that of a finitely generated group $%
G $ and two finitely generated subgroups $\Lambda $ and $\Sigma ,$ not
necessarily cyclic, such that the number of ends of each of the pairs $%
(G,\Lambda )$ and $(G,\Sigma )$ is more than one. Note that any infinite
cyclic subgroup $\Lambda $ of $\pi _{1}(F)$ satisfies $e(\pi _{1}(F),\Lambda
)=2.$ This is because $F$ is closed and orientable so that the cover of $F$
with fundamental group $\Lambda $ is an open annulus which has two ends. In
order to handle the general situation, we will need the concept of an almost
invariant set, which is closely related to the theory of ends. We should
note that Kropholler and Roller \cite{Kropholler} introduced an intersection
cohomology class in the special case of $PD(n-1)$--subgroups of $PDn$--groups.
Their ideas are closely related to ours, and we will discuss the connections
at the start of the last section of this paper. Finally, we should point out
that since Rips asked the above question about symmetry of intersection
numbers, Dunwoody and Sageev \cite{Dunwoody-Sageev} have given a proof of
the existence of a JSJ decomposition for any finitely presented group which
is very much simpler and more elementary than that of Rips and Sela.

The preceding discussion is a little misleading, as the intersection numbers
which we define are not determined simply by a choice of subgroups. In fact,
we define intersection numbers for almost invariant sets. A special case
occurs when one has a group $G$ and subgroups $\Lambda $ and $\Sigma $ such
that $G$ splits over each, as a splitting of $G$ has a well defined almost
invariant set associated. This is discussed in section \ref{intersection
numbers}. Thus we can define the intersection number of two splittings of $G.
$ In the case of cyclic subgroups of surface groups corresponding to simple
closed curves, these curves determine splittings of the surface group over
each cyclic subgroup, and the intersection number we define for these
splittings is the same as the topological intersection number of the curves.

In the first section of this paper, we discuss in more detail intersection
numbers of closed curves on surfaces. In the second section we introduce the
concept of an almost invariant set and prove the symmetry results advertised
in the title. In the third section, we discuss the interpretation of
intersection numbers when they are defined, and how our ideas are connected
with those of Kropholler and Roller.

\rk{Acknowledgments}This paper was written while the author was visiting the
Mathematical Sciences Research Institute in Berkeley in 1996/7. Research at
MSRI is supported in part by NSF grant DMS-9022140. He is also grateful for
the partial support provided by NSF grants DMS-9306240 and DMS-9626537.

\section{The symmetry for surface groups}

\label{symmetry for surface groups}In this section, we will discuss further
the special case of two essential closed curves $L$ and $S$ on a compact
surface $F.$ This will serve to motivate the definitions in the following
section, and also show that the results of that section do indeed answer the
question of Rips. It is not necessary to assume that $F$ is closed or
orientable, but we do need to assume that $L$ and $S$ are two-sided on $F.$
As described in the introduction in the case of simple curves, one defines
their intersection number to be the least number of intersection points
obtainable by homotoping $L$ and $S$ transverse to each other, where the
count is to be made without any signs attached to the intersection points.
(One should also insist that $L$ and $S$ be in general position, in order to
make the count correctly.) Of course, this number is symmetric, ie the
roles of $L$ and $S$ are interchangeable. We will show in section \ref
{intersection numbers} that one can define these intersection numbers in an
algebraically natural way. There is also an idea of self-intersection number
for a curve on a surface and we will discuss a corresponding algebraic idea.

For the next discussion, we will restrict our attention to the case when $L$
and $S$ are simple and introduce the algebraic approach to defining
intersection numbers taken by Rips and Sela in \cite{Rips-Sela}. Let $G$
denote $\pi _{1}(F).$ Suppose that $L$ and $S$ cannot be made disjoint and
choose a basepoint on $L\cap S$. Suppose that $L$ represents the element $%
\lambda $ of $G.$ This element $\lambda $ cannot be trivial, nor can $L$ be
parallel to a boundary component of $F,$ because of our assumption that $L$
and $S$ cannot be made disjoint. Thus $L$ induces a splitting of $G$ over
the infinite cyclic subgroup $\Lambda $ of $G$ which is generated by $%
\lambda .$ Let $\sigma $ denote the element of $G$ represented by $S.$
Define $d(\sigma ,\lambda )$ to be the length of $\sigma $ when written as a
word in cyclically reduced form in the splitting of $G$ determined by $L.$
Similarly, define $d(\lambda ,\sigma )$ to be the length of $\lambda $ when
written as a word in cyclically reduced form in the splitting of $G$
determined by $S.$ For convenience, suppose also that $L$ and $S$ are
separating. Then each of these numbers is equal to the intersection number
of $L$ and $S$ described above and therefore $d(\lambda ,\sigma )=d(\sigma
,\lambda ).$ What is interesting is that this symmetry is not obvious from
the purely algebraic point of view, but it is obvious topologically because
the intersection of two sets is symmetric.

In the above discussion, we restricted attention to simple closed curves on
a surface $F,$ because the algebraic analogue is clear. If $F$ is closed,
then not only does a simple closed curve on $F$ determine a splitting of $%
\pi _{1}(F)$ over the infinite cyclic subgroup carried by the curve, but any
splitting of $\pi _{1}(F)$ over an infinite cyclic subgroup is induced in
this way by some simple closed curve on $F.$ Hence the algebraic situation
described above exactly corresponds to the topological situation when $F$ is
closed.

Now we continue with further discussion of the intersection number of two
closed curves $L$ and $S$ which need not be simple. As in \cite
{FHScurvesonsurfaces}, it will be convenient to assume that $L$ and $S$ are
shortest closed geodesics in some Riemannian metric on $F$ so that they
automatically intersect minimally. Instead of defining the intersection
number of $L$ and $S$ in the ``obvious'' way, we will interpret our
intersection numbers in suitable covers of $F,$ exactly as in \cite
{FHScurvesonsurfaces} and \cite{FHSleastareasurfaces}. Let $F_{\Lambda }$
denote the cover of $F$ with fundamental group equal to $\Lambda .$ Then $L$
lifts to $F_{\Lambda }$ and we denote its lift by $L$ again. Let $l$ denote
the pre-image of this lift in the universal cover $\widetilde{F}$ of $F.$
The full pre-image of $L$ in $\widetilde{F}$ consists of disjoint lines
which we call $L$--lines, which are all translates of $l$ by the action of $%
G. $ Similarly, we define $F_{\Sigma },$ the line $s$ and $S$--lines in $%
\widetilde{F}.$ Now we consider the images of the $L$--lines in $F_{\Sigma }$%
. Each $L$--line has image in $F_{\Sigma }$ which is a possibly singular line
or circle. Then we define $d(L,S)$ to be the number of images of $L$--lines
in $F_{\Sigma }$ which meet $S.$ Similarly, we define $d(S,L)$ to be the
number of images of $S$--lines in $F_{\Lambda }$ which meet $L.$ It is shown
in \cite{FHScurvesonsurfaces}, using the assumption that $L$ and $S$ are
shortest closed geodesics, that each $L$--line in $F_{\Sigma }$ crosses $S$
at most once, and similarly for $S$--lines in $F_{\Lambda }.$ It follows that 
$d(L,S)$ and $d(S,L)$ are each equal to the number of points of $L\cap S,$
and so they are equal to each other. (This assumes that $L$ and $S$ are in
general position.)

Here is an argument which shows that $d(L,S)$ and $d(S,L)$ are equal without
reference to the situation in the surface $F.$ Recall that the $L$--lines are
translates of $l$ by elements of $G.$ Of course, there is not a unique
element of $G$ which sends $l$ to a given $L$--line. In fact, the $L$--lines
are in natural bijective correspondence with the cosets $g\Lambda $ of $%
\Lambda $ in $G.$ (Our groups act on the left on covering spaces.) The
images of the $L$--lines in $F_{\Sigma }$ are in natural bijective
correspondence with the double cosets $\Sigma g\Lambda ,$ and $d(L,S)$
counts the number of these double cosets such that the line $gl$ crosses $s.$
Similarly, $d(S,L)$ counts the number of the double cosets $\Lambda h\Sigma $
such that the line $hs$ crosses $l.$ Note that it is trivial that $gl$
crosses $s$ if and only if $l$ crosses $g{\ }^{-1}s.$ Now we use the
bijection from $G$ to itself given by sending each element to its inverse.
This induces a bijection between the set of all double cosets $\Sigma
g\Lambda $ and the set of all double cosets $\Lambda h\Sigma $ by sending $%
\Sigma g\Lambda $ to $\Lambda g^{-1}\Sigma .$ It follows that it also
induces a bijection between those double cosets $\Sigma g\Lambda $ such that 
$gl$ crosses $s$ and those double cosets $\Lambda h\Sigma $ such that $hs$
crosses $l,$ which shows that $d(L,S)$ equals $d(S,L)$ as required.

This argument has more point when one applies it to a more complicated
situation than that of curves on surfaces. In \cite{FHSleastareasurfaces},
we considered least area maps of surfaces into a 3--manifold. The
intersection number which we used there was defined in essentially the same
way but it had no obvious topological interpretation such as the number of
double curves of intersection. We proved that our intersection numbers were
symmetric by the above double coset argument, in \cite{FHSleastareasurfaces}
just before Theorem 6.3.
\newpage
\section{Intersection Numbers in General}

\label{intersection numbers}In order to handle the general case, we will
need the idea of an almost invariant set. This idea was introduced by Cohen
in \cite{Cohen} and was first used in the relative context by Houghton in 
\cite{Houghton}. We will introduce this idea and explain its connection with
the foregoing.

Let $E$ and $F$ be sets. We say that $E$ and $F$ are almost equal, and write 
$E\stackrel{a}{=}F,$ if the symmetric difference $(E-F)\cup (F-E)$ is
finite. If $E$ is contained in some set $W$ on which a group $G$ acts on the
right, we say that $E$ is almost invariant if $Eg\stackrel{a}{=}E,$ for all $%
g$ in $G.$ An almost invariant subset $E$ of $W$ will be called non-trivial
if it is infinite and has infinite complement. The connection of this idea
with the theory of ends of groups is via the Cayley graph $\Gamma $ of $G$
with respect to some finite generating set of $G.$ (Note that in this paper
groups act on the left on covering spaces and, in particular, $G$ acts on
its Cayley graph on the left.) Using $\mathbb{Z}_{2}$ as coefficients, we
can identify 0--cochains and 1--cochains on $\Gamma $ with sets of vertices or
edges. A subset $E$ of $G$ represents a set of vertices of $\Gamma $ which
we also denote by $E,$ and it is a beautiful fact, due to Cohen \cite{Cohen}%
, that $E$ is an almost invariant subset of $G$ if and only if $\delta E$ is
finite, where $\delta $ is the coboundary operator. If $H$ is a subgroup of $%
G,$ we let $H\backslash G$ denote the set of cosets $Hg$ of $H$ in $G,$ ie
the quotient of $G$ by the left action of $H.$ Of course, $G$ will no longer
act on the left on this quotient, but it will still act on the right. Thus
we have the idea of an almost invariant subset of $H\backslash G.$

Now we again consider the situation of simple closed curves $L$ and $S$ on a
compact surface $F$ and let $\widetilde{F}$ denote the universal cover of $%
F. $ Pick a generating set for $G$ which can be represented by a bouquet of
circles embedded in $F.$ We will assume that the wedge point of the bouquet
does not lie on $L$ or $S.$ The pre-image of this bouquet in $\widetilde{F}$
will be a copy of the Cayley graph $\Gamma $ of $G$ with respect to the
chosen generating set. The pre-image in $F_{\Lambda }$ of the bouquet will
be a copy of the graph $\Lambda \backslash \Gamma ,$ the quotient of $\Gamma 
$ by the action of $\Lambda $ on the left. Consider the closed curve $L$ on $%
F_{\Lambda }.$ Let $D$ denote the set of all vertices of $\Lambda \backslash
\Gamma $ which lie on one side of $L.$ Then $D$ has finite coboundary, as $%
\delta D$ equals exactly the edges of $\Lambda \backslash \Gamma $ which
cross $L.$ Hence $D$ is an almost invariant subset of $\Lambda \backslash G.$
Let $X$ denote the pre-image of $D$ in $\Gamma ,$ so that $X$ equals the set
of vertices of $\Gamma $ which lie on one side of the line $l.$ There is an
algebraic description of $X$ in terms of canonical forms for elements of $G$
as follows. Suppose that $L$ separates $F,$ so that $G=A*_{\Lambda }B.$ Also
suppose that $L$ and $D$ are chosen so that all the vertices of $\Gamma $
labelled with an element of $\Lambda $ do not lie in $X.$ Pick right
transversals $T$ and $T^{\prime }$ for $\Lambda $ in $A$ and $B$
respectively, both of which contain the identity $e$ of $G.$ (A right
transversal of $\Lambda $ in $A$ consists of a choice of coset
representative for each coset $a\Lambda .)$ Each element of $G$ can be
expressed uniquely in the form $a_{1}b_{1}\ldots a_{n}b_{n}\lambda $, where $%
n\geq 1,$ $\lambda $ lies in $\Lambda ,$ each $a_{i}$ lies in $T-\{e\}$
except that $a_{1}$ may be trivial, and each $b_{i}$ lies in $T^{\prime
}-\{e\}$ except that $b_{n}$ may be trivial. Then $X$ consists of those
elements for which $a_{1}$ is non-trivial. If $\Lambda $ is non-separating
in $F,$ there is a similar description for $X.$ See Theorem 1.7 of \cite
{Scott-Wall} for details. Similarly, we can define a set $E$ in $F_{\Sigma }$
and its pre-image $Y$ in $\widetilde{F}$ which equals the set of vertices of 
$\Gamma $ which lie on one side of the line $s.$ Now finally the connection
between the earlier arguments and almost invariant sets can be given. For we
can decide whether the lines $l$ and $s$ cross by considering instead the
sets $X$ and $Y.$ The lines $l$ and $s$ together divide $G$ into the four
sets $X\cap Y,X^{*}\cap Y,X\cap Y^{*}$ and $X^{*}\cap Y^{*},$ where $X{\ }%
^{*}$ denotes $G-X,$ and $l$ crosses $s$ if and only if each of these four
sets projects to an infinite subset of $\Sigma \backslash G.$ Equally, $s$
crosses $l$ if and only if each of these four sets projects to an infinite
subset of $\Lambda \backslash G.$ As we know that $l$ crosses $s$ if and
only if $s$ crosses $l,$ it follows that these conditions are equivalent. We
will show that this symmetry holds in a far more general context.

Note that in the preceding example the subset $X$ of $G$ is $\Lambda $%
--invariant under the left action of $\Lambda $ on $G,$ ie $\lambda X=X,$
for all $\lambda $ in $\Lambda .$

For the most general version of this symmetry result, we can consider any
finitely generated group $G.$ Note that the subgroups of $G$ which we
consider need not be finitely generated.

\begin{definition} 
If $G$ is a finitely generated group and $H$ is a subgroup, then a
subset $X$ of $G$ is {\sl $H$--almost invariant} if $X$ is invariant
under the left action of $ H,$ and simultaneously the quotient set
$H\backslash X$ is almost invariant under the right action of $G.$ In
addition, $X$ is a {\sl non-trivial} $H$--almost invariant subset of $G$ if
$H\backslash X$ and $H\backslash X^{\ast }$ are both infinite.
\end{definition}

Note that if $X$ is a non-trivial $H$--almost invariant subset of $G,$ then $%
e(G,H)$ is at least $2,$ as $H\backslash X$ is a non-trivial almost
invariant subset of $H\backslash G.$

\begin{definition} 
Let $X$ be a $\Lambda $--almost invariant subset of $G$ and let $Y$ be
a $\Sigma $--almost invariant subset of $G.$ We will say that $X$
{\sl crosses} $Y$ if each of the four sets $X\cap Y,X^{\ast }\cap
Y,X\cap Y^{\ast }$ and $X^{\ast }\cap Y^{\ast }$ projects to an
infinite subset of $\Sigma \backslash G.$
\end{definition}

Note that it is obvious that if $Y$ is trivial, then $X$ cannot cross $Y.$
Our first and most basic symmetry result is the following. This is
essentially proved in Lemma 2.3 of \cite{ScottTorusTheorem}, but the context
there is less general.

\begin{lemma} 
\label{crossing is symmetric}If $G$ is a finitely generated group with
subgroups $\Lambda $ and $\Sigma ,$ and $X$ is a non-trivial $\Lambda $%
--almost invariant subset of $G$ and $Y$ is a non-trivial $\Sigma $--almost
invariant subset of $G,$ then $X$ crosses $Y$ if and only if $Y$ crosses $X.$
\end{lemma}

\begin{remark} 
If $X$ and $Y$ are both trivial, then neither can cross the other, so the
above symmetry result is clear. However, this symmetry result fails if only
one of $X$ or $Y$ is trivial. Here is a simple example. Let $\Lambda $ and $%
\Sigma $ denote infinite cyclic groups with generators $\lambda $ and $%
\sigma $ respectively, and let $G$ denote the group $\Lambda \times \Sigma .$
We identify $G$ with the set of integer points in the plane. Let $%
X=\{(m,n)\in G:n>0\},$ and let $Y=\{(m,n)\in G:m=0\}.$ Then $X$ is a
non-trivial $\Lambda $--almost invariant subset of $G$ and $Y$ is a trivial $%
\Sigma $--almost invariant subset of $G.$ One can easily check that $Y$
crosses $X,$ although $X$ cannot cross $Y$ as $Y$ is trivial.
\end{remark}

\begin{proof}
Suppose that $X$ does not cross $Y.$ By replacing one or both of $X$ and $Y$
by its complement if needed, we can assume that $X\cap Y$ projects to a
finite subset of $\Sigma \backslash G.$ The fact that $Y$ is non-trivial
implies that $\Sigma \backslash Y$ is an infinite subset of $\Sigma
\backslash G,$ so there is a point $z$ in $\Sigma \backslash Y$ which is not
in the image of $X\cap Y.$ Now we need to use some choice of generators for $%
G$ and consider the corresponding Cayley graph $\Gamma $ of $G.$ The
vertices of $\Gamma $ are identified with $G$ and the action of $G$ on
itself on the left extends to an action on $\Gamma .$ We consider $z$ and
the image of $X\cap Y$ in the quotient graph $\Sigma \backslash \Gamma .$ As 
$X\cap Y$ has finite image, there is a number $d$ such that each point of
its image can be joined to $z$ by a path of length at most $d.$ As the
projection of $\Gamma $ to $\Sigma \backslash \Gamma $ is a covering map, it
follows that each point of $X\cap Y$ can be joined to some point lying above 
$z$ by a path of length at most $d.$ As any point above $z$ lies in $X^{*},$
it follows that each point of $X\cap Y$ can be joined to some point of $%
X^{*} $ by a path of length at most $d.$ Hence each point of $X\cap Y$ lies
at most distance $d$ from $\delta X.$ Thus the image of $X\cap Y$ in $%
\Lambda \backslash \Gamma $ lies within the $d$--neighbourhood of the compact
set $\delta (\Lambda \backslash X),$ and so must itself be finite. It
follows that $Y$ does not cross $X,$ which completes the proof of the
symmetry result.
\end{proof}

At the start of this section, we explained how to connect the topological
intersection of simple closed curves on a surface with crossing of sets. One
can construct many other interesting examples in much the same way.

\begin{example} 
As before, let $F$ denote a closed surface with fundamental group $G,$ and
let $\widetilde{F}$ denote the universal cover of $F.$ Pick a generating set
of $G$ which can be represented by a bouquet of circles embedded in $F,$ so
that $\widetilde{F}$ contains a copy of the Cayley graph $\Gamma $ of $G$
with respect to the chosen generators. Let $F_{1}$ denote a cover of $F$
which is homeomorphic to a four punctured torus and let $\Lambda $ denote
its fundamental group. For example, if $F$ is the closed orientable surface
of genus three, we can consider a compact subsurface $F^{\prime }$ of $F$
which is homeomorphic to a torus with four open discs removed, and take the
cover $F_{1}$of $F$ such that $\pi _{1}(F_{1})=\pi _{1}(F^{\prime }).$ For
notational convenience, we identify $F_{1}$ with $S^{1}\times S^{1}$ with
the four points $(1,1),(1,i),(1,-1)$ and $(1,-i)$ removed. Now we choose
1--dimensional submanifolds of $F_{1}$ each consisting of two circles and
each separating $F_{1}$ into two pieces. Let $L$ denote $S^{1}\times
\{e^{\pi i/4},e^{5\pi i/4}\}$ and let $S$ denote $S^{1}\times \{e^{3\pi
i/4},e^{7\pi i/4}\}.$ As before, we let $D$ denote all the vertices of the
graph $\Lambda \backslash \Gamma $ in $F_{1}$ which lie on one side of $L,$
and let $E$ denote all the vertices of the graph $\Lambda \backslash \Gamma $
in $F_{1}$ which lie on one side of $S.$ Let $X$ and $Y$ denote the
pre-images of $D$ and $E$ in $G.$ Now $D$ is an almost invariant subset of $%
\Lambda \backslash G,$ as $\delta D$ equals exactly the edges of $\Lambda
\backslash \Gamma $ which cross $L,$ and $E$ is almost invariant for similar
reasons. Hence $X$ and $Y$ are each $\Lambda $--almost invariant subsets of $%
G.$ Clearly $X$ and $Y$ cross. An important feature of this example is that
although $X$ and $Y$ cross, the boundaries $L$ and $S$ of the corresponding
surfaces in $F_{1}$ are disjoint. This is quite different from the example
with which we introduced almost invariant sets, but this is a much more
typical situation.
\end{example}

\begin{definition} 
Let $\Lambda $ and $\Sigma $ be subgroups of a finitely generated group $G.$
Let $D$ denote a non-trivial almost invariant subset of $\Lambda \backslash
G,$ let $E$ denote a non-trivial almost invariant subset of $\Sigma
\backslash G$ and let $X$ and $Y$ denote the pre-images in $G$ of $D$ and $E$
respectively. We define $i(D,E)$ to equal the number of double cosets $%
\Sigma g\Lambda $ such that $gX$ crosses $Y.$
\end{definition}

For this definition to be interesting, we need to show that $i(D,E)$ is
finite, which is not obvious from the definition in this general situation.
In fact, it may well be false if one does not assume that the groups $%
\Lambda $ and $\Sigma $ are finitely generated, although we have no
examples. From now on, we will assume that $\Lambda $ and $\Sigma $ are
finitely generated.

\begin{lemma} 
Let $\Lambda $ and $\Sigma $ be finitely generated subgroups of a finitely
generated group $G.$ Let $D$ denote a non-trivial almost invariant subset of 
$\Lambda \backslash G,$ and let $E$ denote a non-trivial almost invariant
subset of $\Sigma \backslash G.$ Then $i(D,E)$ is finite.
\end{lemma}

\begin{proof}
This is again proved by using the Cayley graph, so it appears to depend on
the fact that $G$ is finitely generated. However, we have no examples where $%
i(D,E)$ is not finite when $G$ is not finitely generated. The proof we give
is essentially contained in that of Lemmas 4.3 and 4.4 of \cite
{Scottendsofpairs}. Start by considering the finite graph $\delta D$ in $%
\Lambda \backslash \Gamma .$ As $\Lambda $ is finitely generated, we can add
edges and vertices to $\delta D$ to obtain a finite connected subgraph $%
\delta _{1}D$ of $\Lambda \backslash \Gamma $ which contains $\delta D$ and
has the property that its inclusion in $\Lambda \backslash \Gamma $ induces
a surjection of its fundamental group to $\Lambda .$ Thus the pre-image of $%
\delta _{1}D$ in $\Gamma $ is a connected graph which we denote by $\delta
_{1}X.$ Similarly, we obtain a finite connected graph $\delta _{1}E$ of $%
\Sigma \backslash \Gamma $ which contains $\delta E$ and has connected
pre-image $\delta _{1}Y$ in $\Gamma .$ As usual, we will denote the
pre-images of $D$ and $E$ in $G$ by $X$ and $Y$ respectively.

Next we claim that if $gX$ crosses $Y$ then $g(\delta _{1}X)$ intersects $%
\delta _{1}Y.$ (The converse need not be true.) Suppose that $g(\delta
_{1}X) $ and $\delta _{1}Y$ are disjoint. Then $g(\delta _{1}X)$ cannot meet 
$\delta Y.$ As $g(\delta _{1}X)$ is connected, it must lie in $Y$ or $Y^{*}.$
It follows that $g(\delta X)$ lies in $Y$ or $Y^{*},$ so that one of the
four sets $X\cap Y,X^{*}\cap Y,X\cap Y^{*}$ and $X^{*}\cap Y^{*}$ must be
empty, which implies that $gX$ does not cross $Y.$

Now we can show that $i(D,E)$ must be finite. Recall that $i(D,E)$ is
defined to be the number of double cosets $\Sigma g\Lambda $ such that $gX$
crosses $Y.$ The preceding paragraph implies that $i(D,E)$ is bounded above
by the number of double cosets $\Sigma g\Lambda $ such that $g(\delta _{1}X)$
meets $\delta _{1}Y.$ Let $P$ and $Q$ be finite subgraphs of $\delta _{1}X$
and $\delta _{1}Y$ which project onto $\delta _{1}D$ and $\delta _{1}E$
respectively. If $g(\delta _{1}X)$ meets $\delta _{1}Y,$ then there exist
elements $\lambda $ of $\Lambda $ and $\sigma $ of $\Sigma $ such that $%
g(\lambda P)$ meets $\sigma Q.$ Thus $\sigma ^{-1}g\lambda P$ meets $Q.$ Now
there are only finitely many elements of $G$ which can translate $P$ to meet 
$Q,$ and it follows that $i(D,E)$ is bounded above by this number.
\end{proof}

We have just shown that, as in the preceding section, the intersection
numbers we have defined are symmetric, but we will need a little more
information.

\begin{lemma} 
Let $G$ be a finitely generated group with subgroups $\Lambda $ and $\Sigma ,
$ let $D$ denote a non-trivial almost invariant subset of $\Lambda
\backslash G,$ and let $E$ denote a non-trivial almost invariant subset of $%
\Sigma \backslash G.$ Then the following statements hold:

\begin{items} 
\item[\bf1\rm)]  $i(D,E)=i(E,D),$

\item[\bf2\rm)]  $i(D,E)=i(D^{\ast },E)=i(D,E^{\ast })=i(D^{\ast },E^{\ast }),$

\item[\bf3\rm)]  if $D^{\prime }$ is almost equal to $D$ and $E^{\prime }$ is almost
equal to $E,$ and $X,X^{\prime }$ and $Y,Y^{\prime }$ denote their
pre-images in $G,$ then $X$ crosses $Y$ if and only if $X^{\prime }$ crosses 
$Y^{\prime },$ so that $i(D,E)=i(D^{\prime },E^{\prime }).$
\end{items}
\end{lemma}

\begin{proof}
The first part is proved by using the bijection from $G$ to itself given by
sending each element to its inverse. This induces a bijection between all
double cosets $\Sigma g\Lambda $ and $\Lambda h\Sigma $ by sending $\Sigma
g\Lambda $ to $\Lambda g^{-1}\Sigma ,$ and it further induces a bijection
between those double cosets $\Sigma g\Lambda $ such that $gX$ crosses $Y$
and those double cosets $\Lambda h\Sigma $ such that $hY$ crosses $X.$

The second part is clear from the definitions.

For the third part, we note that, as $E$ and $E^{\prime }$ are almost equal,
so are their complements in $\Sigma \backslash G,$ and it follows that $X$
crosses $Y$ if and only if it crosses $Y^{\prime }.$ Hence the symmetry
proved in Lemma \ref{crossing is symmetric}, shows that $Y$ crosses $X$ if
and only $Y^{\prime }$ crosses $X.$ Now the same argument reversing the
roles of $D$ and $E$ yields the required result.
\end{proof}

At this point, we have defined in a natural way a number which can
reasonably be called the intersection number of $D$ and $E,$ but have not
yet defined an intersection number for subgroups of $G.$ First note that if $%
e(G,\Lambda )$ is equal to $2,$ then all choices of non-trivial almost
invariant sets in $\Lambda \backslash G$ are almost equal or almost
complementary. Let $D$ denote some choice here. Suppose that $e(G,\Sigma )$
is also equal to $2,$ and let $E$ denote a non-trivial almost invariant
subset of $\Sigma \backslash G.$ The third part of the preceding lemma
implies that $i(D,E)$ is independent of the choices of $D$ and $E$ and so
depends only on the subgroups $\Lambda $ and $\Sigma .$ This is then the
definition of the intersection number $i(\Lambda ,\Sigma ).$ In the special
case when $G\;$is the fundamental group of a closed orientable surface and $%
\Lambda $ and $\Sigma $ are cyclic subgroups of $G,$ it is automatic that $%
e(G,\Lambda )$ and $e(G,\Sigma )$ are each equal to $2.$ The discussion of
the previous section clearly shows that this definition coincides with the
topological definition of intersection number of loops representing
generators of these subgroups, whether or not those loops are simple. Note
that one can also define the self-intersection number of an almost invariant
subset $D$ of $\Lambda \backslash G$ to be $i(D,D),$ and hence can define
the self-intersection number of a subgroup $\Lambda $ of $G$ such that $%
e(G,\Lambda )=2.$ Again this idea generalises the topological idea of
self-intersection number of a loop on a surface.

If one considers subgroups $\Lambda $ and $\Sigma $ such that $e(G,\Lambda )$
or $e(G,\Sigma )$ is greater than $2,$ there are possibly different ideas
for their intersection number depending on which almost invariant sets we
pick. (It is tempting to simply define $i(\Lambda ,\Sigma )$ to be the
minimum possible value for $i(D,E),$ where $D$ is a non-trivial $\Lambda $%
--almost invariant subset of $G$ and $E$ is a non-trivial $\Sigma $--almost
invariant subset of $G.$ But this does not seem to be the ``right''
definition.) However, there is a natural way to choose these almost
invariant sets if we are given splittings of $G$ over $\Lambda $ and $\Sigma
.$ As discussed in the previous section in the case of surface groups, the
standard way to do this when $G=A\ast _{\Lambda }B$ is in terms of canonical
forms for elements of $G$ as follows. Pick right transversals $T$ and $%
T^{\prime }$ for $\Lambda $ in $A$ and $B$ respectively, both of which
contain the identity $e$ of $G.$ Then each element can be expressed uniquely
in the form $a_{1}b_{1}\ldots a_{n}b_{n}\lambda $, where $n\geq 1,$ $\lambda 
$ lies in $\Lambda ,$ each $a_{i}$ lies in $T-\{e\}$ except that $a_{1}$ may
be trivial, and each $b_{i}$ lies in $T^{\prime }-\{e\}$ except that $b_{n}$
may be trivial. Let $X$ denote the subset of $G$ consisting of elements for
which $a_{1}$ is non-trivial, and let $D$ denote $\Lambda \backslash X.$ It
is easy to check directly that $X$ is $\Lambda $--almost invariant. One must
check that $\lambda X=X,$ for all $\lambda $ in $\Lambda $ and that $Dg%
\stackrel{a}{=}D,$ for all $g$ in $G.$ The first equation is trivial, and
the second is easily checked when $g$ lies in $A$ or $B,$ which implies that
it holds for all $g$ in $G.$ Note also that the definition of $X$ is
independent of the choices of transversals of $\Lambda $ in $A$ and $B.$
Then $D$ is the almost invariant set determined by the given splitting of $%
G. $ This definition seems asymmetric, but if instead we consider the $%
\Lambda $--almost invariant subset of $G$ consisting of elements whose
canonical form begins with a non-trivial element of $B,$ we will obtain an
almost invariant subset of $\Lambda \backslash G$ which is almost equal to
the complement of $D.$ There is a similar description of $D$ when $G=A\ast
_{\Lambda }.$ For details see Theorem 1.7 of \cite{Scott-Wall}. The
connection between $D$ and the given splitting of $G$ can be seen in several
ways. From the topologists' point of view, one sees this as described
earlier for surface groups. From the point of view of groups acting on
trees, there is also a very natural description. One identifies a splitting
of $G$ with an action of $G$ on a tree $T$ without inversions, such that the
quotient $G\backslash T$ has a single edge. Let $e$ denote the edge of $T$
with stabiliser $\Lambda ,$ let $v$ denote the vertex of $e$ with stabiliser 
$A,$ and let $E$ denote the component of $T-\{e\}$ which contains $v.$ Then
we can define $X=\{g\in G:ge\subset E\}.$ It is easy to check directly that
this set is the same as the set $X$ defined above using canonical forms.

In the preceding paragraph, we showed how to obtain a well defined
intersection number of given splittings over $\Lambda $ and $\Sigma .$ An
important point to notice is that this intersection number is not determined
by the subgroups $\Lambda $ and $\Sigma $ of $G$ only. It depends on the
given splittings. In the case when $G$ is a surface group, this is
irrelevant as there can be at most one splitting of a surface group over a
given infinite cyclic subgroup. But in general, a group $G$ with subgroup $%
\Lambda $ can have many different splittings over $\Lambda .$

\begin{example} 
Here is a simple example to show that intersection numbers depend on
splittings, not just on subgroups. First we note that the self-intersection
number of any splitting is zero. Now construct a group $G$ by amalgamating
four groups $G_{1},$ $G_{2},G_{3}$ and $G_{4}$ along a common subgroup $%
\Lambda .$ Thus $G$ can be expressed as $G_{12}\ast _{\Lambda }G_{34},$
where $G_{ij}$ is the subgroup of $G$ generated by $G_{i}$ and $G_{j},$ but
it can also be expressed as $G_{13}\ast _{\Lambda }G_{24}$ or $G_{14}\ast
_{\Lambda }G_{23}.$ The intersection number of any distinct pair of these
splittings of $G$ is non-zero, but all the splittings being considered are
splittings over the same group $\Lambda .$
\end{example}

A question which arose in our introduction in connection with the work of
Rips and Sela was how the intersection number of two subgroups of a group $G$
alters if one replaces $G$ by a subgroup. In general, nothing can be said,
but in interesting cases one can understand the answer to this question. The
particular case considered by Rips and Sela was of a finitely presented
group $G$ which is expressed as the fundamental group of a graph of groups
with some vertex group being a group $H$ which contains infinite cyclic
subgroups $\Lambda $ and $\Sigma .$ Further $H$ is the fundamental group of
a surface $F$ and $\Lambda $ and $\Sigma $ are carried by simple closed
curves $L$ and $S$ on $F.$ A point deliberately left unclear in our earlier
discussion of their work was that $F$ is not a closed surface. It is a
compact surface with non-empty boundary. The curves $L$ and $S$ are not
homotopic to boundary components and so define splittings of $H.$ The edges
in the graph of groups which are attached to $H$ all carry some subgroup of
the fundamental group of a boundary component of $F.$ This implies that $L$
and $S$ also define splittings of $G.$ It is clear from this picture that
the intersection number of $\Lambda $ and $\Sigma $ should be the same
whether measured in $G$ or in $H,$ as it should equal the intersection
number of the curves $L$ and $S,$ but this needs a little more thought to
make precise. As usual, the first point to make is that we are really
talking about the intersection numbers of the splittings defined by $L$ and $%
S,$ rather than intersection numbers of $\Lambda $ and $\Sigma .$ For the
number of ends $e(H,\Lambda )$ and $e(H,\Sigma )$ are infinite when $F$ is a
surface with boundary. As $G$ is finitely presented, we can attach cells to
the boundary of $F$ to construct a finite complex $K$ with fundamental group 
$G$. Now the identification of the intersection number of the given
splittings of $G$ with the intersection number of $L$ and $S$ proceeds
exactly as at the start of this section, where we showed how to identify the
intersection number of the given splittings of $H$ with the intersection
number of $L$ and $S.$

\section{Interpreting intersection numbers}

\label{more about intersection numbers}It is natural to ask what is the
meaning of the intersection numbers defined in the previous section. The
answer is already clear in the case of a surface group with cyclic
subgroups. In this section, we will give an interpretation of the
intersection number of two splittings of a finitely generated group $G$ over
finitely generated subgroups. We start by discussing the connection with the
work of Kropholler and Roller.

In \cite{Kropholler}, Kropholler and Roller introduced an intersection
cohomology class for $PD(n-1)$--subgroups of a $PDn$--group. The pairs
involved always have two ends, so the work of the previous section defines
an intersection number in this situation. The connection between our
intersection number and their intersection cohomology class is the
following. Recall that if one has subgroups $\Lambda $ and $\Sigma $ of a
finitely generated group $G,$ such that $e(G,\Lambda )$ and $e(G,\Sigma )$
are each equal to $2,$ then one chooses a non-trivial $\Lambda $--almost
invariant subset $X$ of $G$ and a non-trivial $\Sigma $--almost invariant
subset $Y$ of $G$ and defines our intersection number $i(\Lambda ,\Sigma )$
to equal the number of double cosets $\Sigma g\Lambda $ such that $gX$
crosses $Y.$ Their cohomology class encodes the information about which
double cosets have this crossing property. Thus their invariant is much
finer than the intersection number and it is trivial to deduce the
intersection number from their cohomology class.

To interpret the intersection number of two splittings of a group $G,$ we
need to discuss the Subgroup Theorem for amalgamated free products. Let $G$
be a finitely generated group, which splits over finitely generated
subgroups $\Lambda $ and $\Sigma .$ We will write $G=A_{1}*_{\Lambda
}(B_{1}) $ to denote that either $G$ has the HNN structure $A_{1}*_{\Lambda
} $ or $G$ has the structure $A_{1}*_{\Lambda }B_{1}.$ Similarly, we will
write $G=A_{2}*_{\Sigma }(B_{2}).$ The Subgroup Theorem, see \cite
{Scott-Wall} and \cite{Serre} (or \cite{Serre2}) for discussions from the
topological and algebraic points of view, yields a graph of groups structure 
$\Phi _{1}(\Sigma )$ for $\Sigma ,$ with vertex groups lying in conjugates
of $A_{1}$ or $B_{1}$ and edge groups lying in conjugates of $\Lambda .$
Typically this graph will not be finite or even locally finite. However, as $%
\Sigma $ is finitely generated, there is a finite subgraph $\Psi _{1}$ which
still carries $\Sigma .$ If we reverse the roles of $\Lambda $ and $\Sigma ,$
we will obtain a graph of groups structure $\Phi _{2}(\Lambda )$ for $%
\Lambda ,$ with vertex groups lying in conjugates of $A_{2}$ or $B_{2}$ and
edge groups lying in conjugates of $\Sigma ,$ and there is a finite subgraph 
$\Psi _{2}$ which still carries $\Lambda .$ We show below that, in most
cases, the intersection number of $\Lambda $ and $\Sigma $ measures the
minimal possible number of edges of these finite subgraphs. Notice that if
we consider the special case when $G$ is the fundamental group of a closed
surface and $\Lambda $ and $\Sigma $ are infinite cyclic subgroups, this
statement is clear. Now the symmetry of intersection numbers implies the
surprising fact that the minimal number of edges for $\Psi _{1}$ and $\Psi
_{2}$ are the same.

There is an alternative point of view which we will use for our proof. The
splitting $A_{2}*_{\Sigma }(B_{2})$ of $G$ corresponds to an action of $G$
on a tree $T$ such that the quotient $G\backslash T$ has one edge. The edge
stabilisers in this action on $T$ are all conjugate to $\Sigma $ and the
vertex stabilisers are conjugate to $A_{2}$ or $B_{2}$ as appropriate. If
one has a subgroup $\Lambda $ of $G,$ the quotient $\Lambda \backslash T$
will be the graph underlying $\Phi _{2}(\Lambda ).$ There is a $\Lambda $%
--invariant subtree $T^{\prime }$ of $T,$ such that the graph $\Lambda
\backslash T^{\prime }$ is the graph underlying $\Psi _{2}.$ Whichever point
of view you take, it is necessary to connect it with the ideas about almost
invariant sets which we have already discussed. Here is our interpretation
of intersection numbers.

\begin{theorem} 
\label{interpret intersection numbers}Let $G$ be a finitely generated group,
which splits over finitely generated subgroups $\Lambda $ and $\Sigma ,$
such that if $U$ and $V$ are any conjugates of $\Lambda $ and $\Sigma $
respectively, then $U\cap V$ has infinite index in both $U$ and $V.$ Then
the intersection number of the two splittings equals the minimal number of
edges in each of the graphs $\Psi _{1}$ and $\Psi _{2}.$
\end{theorem}

\begin{remark} 
This result is clearly false if the condition on conjugates is omitted. For
example, if $\Lambda =\Sigma ,$ then $\Psi _{1}(\Sigma )$ and $\Psi
_{2}(\Lambda )$ will each consist of a single vertex, but the intersection
number of the two splittings need not be zero.
\end{remark}

The proof will use the following sequence of lemmas.

We start with a general result about minimal $G$--invariant subtrees of a
tree $T$ on which a group $G$ acts. If every element of $G$ fixes each point
of a non-trivial subtree $T^{\prime }$ of $T,$ then any vertex of $T^{\prime
}$ is a minimal $G$--invariant subtree of $T.$ Otherwise, there is a unique
minimal $G$--invariant subtree of $T.$ An orientation of an edge $e$ of $T$
consists of a choice of one vertex as the initial vertex $i(e)$ of $e$ and
the other as the terminal vertex $t(e).$ An oriented path in $T$ consists of
a finite sequence of oriented edges $e_{1},e_{2},\ldots ,e_{k}$ of $T,$ such
that $t(e_{j})=i(e_{j+1}),$ for $1\leq j\leq k-1.$ If we consider two
oriented edges $e$ and $e^{\prime }$ of $T$ we say that they are coherently
oriented if there is an oriented path which begins with one and ends with
the other. Finally, given an edge $e$ of $T$ and an element $g$ of $G$, we
will say that $e$ and $ge$ are coherently oriented if for some (and hence
either) orientation on $e$ and the induced orientation on $ge,$ the edges $e$
and $ge$ are coherently oriented.

\begin{lemma} 
\label{edge lies in minimal subtree}Suppose that a group $G$ acts on a tree $%
T$ without inversions and without fixing a point. Let $T^{\prime }$ denote
the minimal $G$--invariant subtree. Then an edge $e$ of $T$ lies in $%
T^{\prime }$ if and only if there exists an element $g$ of $G$ such that $e$
and $ge$ are distinct and coherently oriented.
\end{lemma}

\begin{proof}
First consider an edge $e$ not lying in $T^{\prime }.$ Orient $e$ so that it
is the first edge of an oriented path $\lambda $ in $T$ which starts with $%
e, $ has no edge in $T^{\prime },$ and ends at a vertex of $T^{\prime }.$
Thus $ge,$ with the induced orientation, is the first edge of an oriented
path $g\lambda $ in $T$ which starts with $ge,$ has no edge in $T^{\prime },$
and ends at a vertex of $T^{\prime }.$ Now the unique path in $T$ which
joins $e$ and $ge$ must consist either of $\lambda $ and $g\lambda $
together with a path in $T^{\prime }$ or of an initial segment of $e$
together with an initial segment of $ge.$ In either case, it follows that $e$
and $ge$ are not coherently oriented.

Now we consider an edge $e$ of $T^{\prime }$ and its image $\overline{e}$ in 
$G\backslash T^{\prime }.$

If $\overline{e}$ is non-separating in $G\backslash T^{\prime },$ let $\mu $
denote an oriented path in $G\backslash T^{\prime }$ which joins the ends of 
$\overline{e}$ and meets $\overline{e}$ only in its endpoints. Then the loop
formed by $\mu \cup \overline{e}$ lifts to an oriented path in $T^{\prime },$
which shows that there is $g$ in $G$ such that $e$ and $ge$ are distinct and
coherently oriented.

If $\overline{e}$ separates $G\backslash T^{\prime },$ we can write the
graph $G\backslash T^{\prime }$ as $\Gamma _{1}\cup \overline{e}\cup \Gamma
_{2},$ where each $\Gamma _{i}$ is connected and meets $\overline{e}$ in one
endpoint only. Now consider the graph of groups structure given by $%
G\backslash T^{\prime }.$ By contracting each $\Gamma _{i}$ to a point, we
obtain an amalgamated free product structure of $G$ as $G_{1}*_{C}G_{2},$
where $C=stab(e)$ and each $G_{i}$ is the fundamental group of the graph of
groups $\Gamma _{i}.$ Let $T_{i}$ denote the tree on which $G_{i}$ acts with
quotient $\Gamma _{i}.$ Then the complement in $T^{\prime }$ of the edge $e$
and its translates consists of disjoint copies of $T_{1}$ and $T_{2}.$ We
identify $T_{i}$ with the copy of $T_{i}$ which meets $e.$ Note that $T_{1}$
and $T_{2}$ are disjoint. Now it is clear that $G_{1}\neq C\neq G_{2}.$ For
if $G_{1}=C,$ then $G=G_{2},$ which implies that $T_{2}$ is a $G$--invariant
subtree of $T^{\prime },$ contradicting the minimality of $T^{\prime }.$ As $%
G_{1}\neq C,$ there is an element $g_{1}$ of $G_{1}$ such that $g_{1}e\neq
e, $ and similarly there is an element $g_{2}$ of $G_{2}$ such that $%
g_{2}e\neq e.$ For each $i,$ there is a path $\lambda _{i}$ in $T_{i}$ which
begins at $e$ and ends at $g_{i}e.$ As $T_{1}$ and $T_{2}$ are disjoint, so
are $\lambda _{1}$ and $\lambda _{2}.$ It follows that of the three edges $%
e,g_{1}e,g_{2}e,$ at least one pair is coherently oriented, which completes
the proof of the lemma.
\end{proof}

The following result is clear.

\begin{lemma} 
\label{connect edges with gE}Suppose that a group $G$ acts on a tree $T$
without inversions and without fixing a point. Let $e$ denote an edge of $T,$
let $E$ denote a component of $T-\{e\}$ and let $g$ denote an element of $G.$
Then $e$ and $ge$ are distinct and coherently oriented if and only if either 
$gE\subsetneqq E$ or $gE^{\ast }\subsetneqq E^{\ast }.$
\end{lemma}

Next we need to connect this with almost invariant sets, although the
following result does not use the almost invariance property.

\begin{lemma} 
\label{connect gE with gY}Suppose that a group $G$ acts on a tree $T$
without inversions and without fixing a point and suppose that the quotient
graph $G\backslash T$ has only one edge. Let $e$ denote an edge of $T,$ let $%
E$ denote a component of $T-\{e\}$ and let $Y=\{k\in G:ke\subset E\}.$ Then
the following statements hold for all elements $g$ of $G:$

\begin{items}
\item[\bf1\rm)]  $gY\subset Y$ if and only if $gE\subset E,$ and $gY^{\ast }\subset
Y^{\ast }$ if and only if $gE^{\ast }\subset E^{\ast }.$

\item[\bf2\rm)]  $gY=Y$ if and only if $gE=E,$ and $gY^{\ast }=Y^{\ast }$ if and only
if $gE^{\ast }=E^{\ast }.$

\item[\bf3\rm)]  $gY\subsetneqq Y$ if and only if $gE\subsetneqq E,$ and $gY^{\ast
}\subsetneqq Y^{\ast }$ if and only if $gE^{\ast }\subsetneqq E^{\ast }.$
\end{items}
\end{lemma}

\begin{proof}
Suppose that $gE\subset E.$ If $k$ lies in $Y,$ then $ke\subset E,$ so that $%
gke\subset gE\subset E.$ Thus $gk$ also lies in $Y.$ It follows that $%
gY\subset Y.$

Conversely, suppose that $gY\subset Y$ and consider an edge $f$ of $E.$ As $%
G\backslash T$ has only one edge, $f=ke$ for some $k$ in $G.$ As $f$ lies in 
$E,$ $k$ lies in $Y,$ and hence $gk$ also lies in $Y$ by our assumption that 
$gY\subset Y.$ Thus $gke\subset E,$ so that $gf\subset E.$ Thus implies that 
$gE\subset E$ as required.

The proof for the second equivalence in part 1 is essentially the same.

The equivalences in part 2 follow by applying part 1 for $g$ and $g^{-1}.$
Now the equivalences in part 3 are clear.
\end{proof}

Next we connect the above inclusions with crossing of sets.

\begin{lemma} 
\label{connects crossing with gY}Suppose that a finitely generated group $G$
splits over a finitely generated subgroup $\Lambda $ with corresponding $%
\Lambda $--almost invariant set $X$ and also splits over a finitely generated
subgroup $\Sigma $ with corresponding $\Sigma $--almost invariant set $Y.$
Suppose further that if $U$ and $V$ are any conjugates of $\Lambda $ and $%
\Sigma $ respectively, then $U\cap V$ has infinite index in $U.$ Then $X$
crosses $Y$ if and only if there is an element $\lambda $ in $\Lambda $ such
that either $\lambda Y\subsetneqq Y$ or $\lambda Y^{\ast }\subsetneqq
Y^{\ast }.$
\end{lemma}

\begin{proof}
We claim that there exists $\lambda _{1}\in \Lambda $ such that either $%
\lambda _{1}Y\subsetneqq Y$ or $\lambda _{1}Y^{\ast }\subsetneqq Y,$ and
there exists $\lambda _{2}\in \Lambda $ such that either $\lambda
_{2}Y\subsetneqq Y^{\ast }$ or $\lambda _{2}Y^{\ast }\subsetneqq Y^{\ast }.$
Assuming this, either $\lambda _{1}Y\subsetneqq Y$ or $\lambda _{2}Y^{\ast
}\subsetneqq Y^{\ast },$ and our proof is complete, or we have $\lambda
_{1}Y^{\ast }\subsetneqq Y$ and $\lambda _{2}Y\subsetneqq Y^{\ast }.$ The
last possibility implies that $\lambda _{2}\lambda _{1}Y^{\ast }\subsetneqq
\lambda _{2}Y\subsetneqq Y^{\ast },$ again completing the proof.

To prove our claim, we pick a finite generating set for $G,$ and consider
the Cayley graph $\Gamma $ of $G$ with respect to this generating set. As $Y$
is a $\Sigma $--almost invariant set associated to a splitting $A_{2}\ast
_{\Sigma }(B_{2})$ of $G$ over $\Sigma ,$ we can choose $\Gamma $ and $Y$ so
that, for every $g$ in $G$, $g\delta Y$ is disjoint from or coincides with $%
\delta Y.$ A simple way to arrange this is to take as generators of $G$ the
union of a set of generators of $\Sigma $ and of $A_{2}$ and $B_{2},$ so
that $\Gamma (G)$ contains a copy of the Cayley graph $\Gamma (\Sigma )$ of $%
\Sigma $ and $\Sigma \backslash \Gamma $ contains $\Sigma \backslash \Gamma
(\Sigma )$ which is a wedge of circles. (Note that this uses the hypothesis
that $\Sigma $ is finitely generated.) Let $v$ denote the wedge point, and
let $E$ denote the collection of vertices of $\Sigma \backslash \Gamma $
which can be joined to $v$ by a path whose interior is disjoint from $v$
such that the last edge is labelled by an element of $A.$ Then clearly $%
\delta E$ consists of exactly those edges of $\Sigma \backslash \Gamma $
which have one end at $v$ and are labelled by an element of $A.$ Further, if
we let $Y$ denote the pre-image of $E$ in $G,$ then, for every $g$ in $G$, $%
g\delta Y$ is disjoint from or coincides with $\delta Y.$

In order to prove that $\lambda _{1}$ exists, we argue as follows. As $%
\Lambda \cap \Sigma $ has infinite index in $\Lambda ,$ and as $\delta X$ is 
$\Lambda $--invariant, it follows that $\delta X$ must contain points which
are arbitrarily far from $\delta Y$ on each side of $\delta Y.$ Recall that $%
\Lambda \backslash X$ is an almost invariant subset of $\Lambda \backslash
G, $ so that it has finite coboundary which equals $\Lambda \backslash
\delta X. $ Hence there is a number $d$ such that any point of $\Lambda
\backslash \delta X$ can be joined to the image of $\delta Y$ in $\Lambda
\backslash \Gamma $ by a path of length at most $d.$ It follows that any
point of $\delta X$ can be joined to $\lambda \delta Y,$ for some $\lambda $
in $\Lambda ,$ by a path in $\Gamma $ of length at most $d.$ Hence there is
a translate of $\delta Y$ which contains points on one side of $\delta Y$
and another translate which contains points on the other side of $\delta Y.$
Hence there are elements $\lambda _{1}$ and $\lambda _{2}$ of $\Lambda $
such that $\lambda _{1}\delta Y$ lies on one side of $\delta Y$ and $\lambda
_{2}\delta Y$ lies on the other. Without loss of generality, we can suppose
that $\lambda _{1}\delta Y$ lies on the side containing $Y$ so that either $%
\lambda _{1}Y\subsetneqq Y$ or $\lambda _{1}Y^{\ast }\subsetneqq Y.$ As $%
\lambda _{2}\delta Y$ lies on the side of $\delta Y$ containing $Y^{\ast },$
either $\lambda _{2}Y\subsetneqq Y^{\ast }$ or $\lambda _{2}Y^{\ast
}\subsetneqq Y^{\ast }.$ This completes the proof of the claim made at the
start of the proof.
\end{proof}

Now we can give the proof of Theorem \ref{interpret intersection numbers}.

\begin{proof}
Recall that $G$ splits over finitely generated subgroups $\Lambda $ and $%
\Sigma $ such that if $U$ and $V$ are any conjugates of $\Lambda $ and $%
\Sigma ,$ then $U\cap V$ has infinite index in both $U$ and $V.$ Also $G$
acts on a tree $T$ so as to induce the given splitting over $\Sigma .$ Let $%
e $ denote an edge of $T$ with stabiliser $\Sigma $ and consider the action
of $\Lambda $ on $T.$ Our hypothesis on conjugates of $\Lambda $ and $\Sigma 
$ implies, in particular, that $\Lambda $ is not contained in any conjugate
of $\Sigma $ so that $\Lambda $ cannot fix an edge of $T.$ Thus there is a
unique minimal $\Lambda $--invariant subtree $T^{\prime }$ of $T.$ Lemma \ref
{edge lies in minimal subtree} shows that an edge $he$ of $T$ lies in $%
T^{\prime }$ if and only if there is $\lambda $ in $\Lambda $ such that $he$
and $\lambda he$ are distinct and coherently oriented. Lemma \ref{connect
edges with gE} shows that this occurs if and only if either $\lambda
hE\subsetneqq hE$ or $\lambda hE^{\ast }\subsetneqq hE^{\ast },$ and Lemma 
\ref{connect gE with gY} shows that this occurs if and only if $\lambda
hY\subsetneqq hY$ or $\lambda hY^{\ast }\subsetneqq hY^{\ast }.$ Finally
Lemma \ref{connects crossing with gY} shows that this occurs if and only if $%
X$ crosses $hY.$ We conclude that an edge $he$ of $T$ lies in $T^{\prime }$
if and only if $X$ crosses $hY.$ Thus the edges of $T$ which lie in the
minimal $\Lambda $--invariant subtree $T^{\prime }$ naturally correspond to
the cosets $h\Sigma $ such that $X$ crosses $hY.$ Hence the number of edges
in $\Psi _{2}(\Lambda )$ equals the number of double cosets $\Lambda h\Sigma 
$ such that $X$ crosses $hY,$ which was defined to be the intersection
number of the given splittings. Similarly, one can show that the
intersection number of the given splittings equals the minimal possible
number of edges in the graph $\Psi _{1}(\Sigma ).$ This completes the proof
of Theorem \ref{interpret intersection numbers}.
\end{proof}

\newpage
\count0=333

\shorttitle{Correction}

\gt\hfill      
\beginpicture
\setcoordinatesystem units <0.33truein, 0.33truein> point at 2.2 0.9
\setplotsymbol ({$\cal G$})
\plotsymbolspacing=9truept
\circulararc 315 degrees from 0 1 center at 0 0
\setplotsymbol ({$\cal T$})
\circulararc 315 degrees from 1 -1 center at 1 0
\endpicture
%
\break
{\small Volume 2 (1998) 333--335\nl
Erratum 1\nl
Published:  7 September 1998}

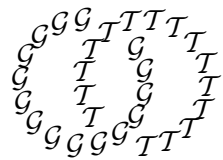
\vglue 0.4truein

\centerline{\Large \bf Correction to ``The Symmetry of Intersection}
\smallskip
\cl{\Large \bf  Numbers in Group Theory''}
\medskip
\cl{\sc Peter Scott}
\bigskip

Theorem 3.1 is false as stated. The error in the argument occurs in the
proof of Lemma 3.6. See below for a counterexample.

Lemma 3.6 asserts that, under suitable hypotheses, $X$ crosses $Y$ if and
only if there is an element $\lambda $ in $\Lambda $ such that either $%
\lambda Y\subsetneqq Y$ or $\lambda Y^{\ast }\subsetneqq Y^{\ast }.$ One of
these implications is correct. If such a $\lambda $ exists, then it is true
that $X$ must cross $Y.$ However, I\ failed to give any argument for this,
and I provide one below. The other implication is false. The mistake is
contained in the second sentence of the second paragraph on page 28. A
simple fix is to amend the statements of Theorem 3.1 and Lemma 3.6 to take
this into account. Thus we need the additional hypothesis for Lemma 3.6 that
if $X$ crosses $Y,$ then $\delta X$ must contain points which are
arbitrarily far from $\delta Y$ on each side of $\delta Y.$ We also need the
additional hypothesis for Theorem 3.1 that if $X$ crosses $gY$ then $\delta
X $ must contain points which are arbitrarily far from $\delta gY$ on each
side of $\delta gY.$ This technical assumption is often but not always
satisfied.

Here is the half of the proof of Lemma 3.6 which was omitted. This asserts
that if there is an element $\lambda $ in $\Lambda $ such that either $%
\lambda Y\subsetneqq Y$ or $\lambda Y^{\ast }\subsetneqq Y^{\ast },$ then $X$
must cross $Y.$ We will assume that $\lambda Y\subsetneqq Y,$ as the
argument in the other case is essentially identical. As $Y$ is associated to
a splitting of $G,$ it is easy to see that the distance of $\lambda
^{n}\delta Y$ from $\delta Y$ must tend to infinity as $n\rightarrow \infty
. $ (For example, if $G=A\ast _{C}B,$ and $Y$ is the set of words in $G$
which begin in $A-C,$ then $\lambda $ must begin in $A-C$ and end in $B-C.)$
Now consider an element $g\in G,$ and let $d$ denote the distance of $g$
from $\delta Y.$ Then $d$ is also the distance of $\lambda ^{n}g$ from $%
\lambda ^{n}\delta Y.$ Hence, for any element $g$ of $G,$ all translates $%
\lambda ^{n}g$ lie in $Y,$ for suitably large $n.$ If we apply these
statements to an edge of $\delta X,$ and recall that $\delta X$ is preserved
by $\lambda ,$ we see that $\delta X$ must contain points which are
arbitrarily far from $\delta Y$ and lie in $Y.$ By applying the same
discussion to $\lambda ^{-1},$ we see that $\delta X$ must also contain
points which are arbitrarily far from $\delta Y$ and lie in $Y^{\ast }.$
Hence $\delta X$ must contain points which are arbitrarily far from $\delta
Y $ on each side of $\delta Y$ as required.

Now we come to the promised counterexample. Let $G$ denote the fundamental
group of the closed orientable surface $M$ of genus two. Let $D$ denote a
simple closed curve on $M$ which separates $M$ into two once-punctured tori $%
S$ and $T$ and let $D^{\prime }$ denote a non-separating simple closed curve
in the interior of $S.$ Let $W$ denote the surface obtained from $S$ by
removing a regular neighbourhood of $D^{\prime }.$ Let $C$ denote a
non-separating simple closed curve on $M$ whose intersection number with $D$
is two, and which is disjoint from $D^{\prime }.$ We will describe two
splittings of $G.$ The first will be the HNN\ splitting over an infinite
cyclic group determined by $C.$ The second will be the amalgamated free
product splitting of $G\;$over $\pi _{1}(W)$ with vertex groups $\pi _{1}(S)$
and $\pi _{1}(W\cup T).$ These two splittings satisfy the hypotheses of
Theorem 3.1. If one considers $\pi _{1}(C)$ as a subgroup of the splitting
over $\pi _{1}(W),$ the minimal graph of groups obtained has no edges,
because $\pi _{1}(C)$ is contained in $\pi _{1}(W\cup T)$ which is a vertex
group. If one considers $\pi _{1}(W)$ as a subgroup of the HNN splitting
determined by $C,$ the minimal graph of groups obtained has at least one
edge because $\pi _{1}(W)$ does not lie in a conjugate of any vertex group.
(The graph in question has exactly one edge, but this fact is not needed
here.) This shows that Theorem 3.1 must fail for this example, because the
numbers of edges in these two graphs are not equal. It also true that Lemma
3.6 fails for this example. Let $X$ and $Y$ be the usual subsets of $G$
associated to the two splittings. I claim that $X$ crosses $Y$ but $\delta X$
does not contain points which are arbitrarily far from $\delta Y$ on each
side of $\delta Y.$ To see this, consider the picture in the cover $M_{C}$
of $M$ whose fundamental group equals $\pi _{1}(C).$ This cover is an open
annulus which contains a lift of $C$ which we will continue to denote by $C.$
As in section 2, we pick a generating set for $G$ which can be represented
by a bouquet of circles embedded in $M,$ so that the pre-image in the
universal cover $\widetilde{M}$ of $M$ of this bouquet is a copy of the
Cayley graph $\Gamma $ of $G,$ and we identify the vertices of this graph
with $G.$ Now let $E$ denote the set of all vertices of $\pi
_{1}(C)\backslash \Gamma $ in $M_{C}$ which lie on one side of $C.$ Then $E$
represents an almost invariant subset of $\pi _{1}(C)\backslash G$ and the
pre-image of $E$ in $\Gamma $ can be taken to be $X.$ Now consider the
picture in the cover $M_{W}$ of $M$ whose fundamental group equals $\pi
_{1}(W).$ This cover consists of a lift of $W,$ which we will continue to
denote by $W$ and open collars attached to the boundary components of $W.$
Let $F$ denote the set of all vertices of $\pi _{1}(W)\backslash \Gamma $
which lie in the union of $W$ together with the collar attached to the
component $D$ of $\partial W.$ Then $F$ represents an almost invariant
subset of $\pi _{1}(W)\backslash G$ and the pre-image of $F$ in $G$ can be
taken to be $Y.$ The pre-image in $\widetilde{M}$ of $C$ is a line whose
image in $M_{W}$ is a properly embedded line meeting $W$ in a compact arc
which projects homeomorphically to $C\cap W.$ Now inspection shows that each
of the four sets $X^{(\ast )}\cap Y^{(\ast )}$ has infinite image in $M_{W}$
so that $X$ crosses $Y$ but $\delta X$ does not contain points which are
arbitrarily far from $\delta Y$ on each side of $\delta Y.$

The new version of Theorem 3.1 described here is, of course, rather
unsatisfactory as the extra hypothesis is technical and it is not clear when
it holds. However, there is a little more which can be said without any
extra work. For it follows from the preceding discussion that the number of
edges in each of the minimal graphs of groups described above is always less
than or equal to the intersection number of the two splittings being
considered.

\end{document}